\documentclass[11pt]{amsart}
\usepackage{amsmath}
\usepackage{amsfonts}
\usepackage{amssymb}
\usepackage{epsfig,latexsym,graphicx}

\newtheorem{Thm}{Theorem}[section]
\newtheorem{thm}[Thm]{Theorem}

\newtheorem{lemma}[Thm]{Lemma}
\newtheorem{proposition}[Thm]{Proposition}

\newtheorem{conjecture}[Thm]{Conjecture}

\newcommand{\ip}[2]{{\langle#1,#2\rangle}}
\newcommand{\ipT}[2]{{\langle#1,#2\rangle}}
\newcommand{\norm}[1]{{\|#1\|}}
\newcommand{\LTFnorm}{{\Lambda_\fc}}

\newcommand{\eps}{\varepsilon}
\newcommand{\qq}{Q}

\def\ldots{\mathinner{\ldotp\ldotp\ldotp}}
\def\ldots{\mathinner{\cdotp\cdotp\cdotp}}

\def \K{\rm {\bf K}}

\def \tr{\mathop{\mathrm{tr}}}

\def \cal{\mathbb}

\def \beq{\begin{eqnarray*}}
\def \eeq{\end{eqnarray*}}

\def \R{\mathbb{R}}

\def \K{\mathbb{K}}
\def \E{\mathbb{\cal E}}

\newcommand{\fc}{\mathcal{F}}
\newcommand{\I}{{\mathbb I}}

\newcommand{\nonlin}{{\alpha_\fc}}
\newcommand{\Ncal}{{\cal N}}
\newcommand{\recon}{{\omega}}
\newcommand{\intOmega}{{\mathring{\Omega}_1}}
\newcommand{\Sc}{\mathcal{S}}
\newcommand{\Ac}{\mathcal{A}}
\newcommand{\Sym}{\mathrm{Sym}}
\newcommand{\sigmabf}{{\boldsymbol\sigma}}
\newcommand{\ignore}[1]{}

\newcommand{\ra}{{\rightarrow}}

\newcommand{\eproof}{\hfill\rule{2.2mm}{3.0mm}}

\newcommand{\Proof}{\noindent {\bf Proof.~~}}

\renewcommand{\eqref}[1]{(\ref{#1})}
\newcommand{\inner}[1]{\langle #1 \rangle}
\newcommand{\shsp}{\hspace{1em}}
\newcommand{\mhsp}{\hspace{2em}}

\newcommand{\rank}{{\rm rank}}

\begin{document}
\baselineskip 14pt

\title{Invertibility and Robustness of Phaseless Reconstruction}
\author[R.~Balan]{Radu Balan}
\address[R.~Balan]{Department of Mathematics, University of Maryland, College Park MD 20742}
\email[R.~Balan]{rvbalan@math.umd.edu}

\author[Y.~Wang]{Yang Wang}
\address[Y.~Wang]{School of Mathematics, Michigan State University, East Lansing, MI 48824}
\email[Y.~Wang]{ywang@math.msu.edu}
\thanks{R.~Balan was supported in part by NSF DMS-1109498. Y.~Wang was supported in part by NSF DMS-1043032, and by AFOSR FA9550-12-1-0455.}


\begin{abstract}
This paper is concerned with the question of reconstructing
a vector in a finite-dimensional real Hilbert space when
only the magnitudes of the coefficients of the vector under a redundant linear map
are known. We analyze various Lipschitz bounds of the nonlinear analysis map and we establish
theoretical performance bounds of any reconstruction algorithm.
We show that robust and stable reconstruction requires
additional redundancy than the critical threshold.
\end{abstract}

\maketitle

\section{Introduction}

This paper is concerned with the question of reconstructing
a vector $x$ in a finite-dimensional real Hilbert space $H$ of dimension $n$ when
only the magnitudes of the coefficients of the vector under a redundant linear map
are known.

Specifically our problem is to reconstruct $x\in H$ up to a global phase factor from the 
magnitudes $\{ |\ip{x}{f_k}|~,~1\leq k\leq m\}$ where 
$\fc=\{f_1,\ldots,f_m\}$ is a frame (complete system) for $H$.

A previous paper \cite{BCE06} described the importance of this problem to signal processing,
in particular to the analysis of speech. Of particular interest is the case when the
coefficients are obtained from a Windowed Fourier Transform (also known as Short-Time Fourier
Transform), or an Undecimated Wavelet Transform (in audio and image signal processing). A
similar problem appears in Quantum Information (QI) literature (see e.g. \cite{HMW11}).
However some important differences are notable: first the unknown objects to be reconstructed
are quantum states (meaning nonnegative symmetric operators of unit trace); secondly,
measurements are performed by taking Hilbert-Schmidt inner products against some (nonnegative)
symmetric operators of rank not necessarily one. In QI language our problem is to reconstruct
rank one nonnegative symmetric operators from measurements against a set of rank one
nonnegative symmetric operators.

 While \cite{BCE06} presents some necessary and
sufficient conditions for reconstruction, the general problem of finding
fast/efficient algorithms is still open.  In \cite{Bal2010} we describe one
solution in the case of STFT coefficients.

For vectors in real Hilbert spaces,  the reconstruction problem
is easily shown to be equivalent to a
combinatorial problem.  In \cite{BCE07} this problem is
further proved to be equivalent to a (nonconvex) optimization problem.

A different approach (which we called the {\em algebraic approach}) was proposed in
\cite{Bal2009}. While it applies to both real and complex cases, noisless and noisy cases, the
approach requires solving a linear system of size exponentially in space dimension.  This
algebraic approach generalizes the approach in \cite{BBCE07} where reconstruction is performed
with complexity $O(n^2)$ (plus computation of the principal eigenvector for a matrix of size
$n$). However this method requires $m=O(n^2)$ frame vectors.

Recently the authors of \cite{CSV12} developed a convex optimization algorithm (a SemiDefinite
Program called {\em PhaseLift}) and proved its ability to perform exact reconstruction in the
absence of noise, as well as its stability under noise conditions. In a separate paper
\cite{CESV12}, the authors further developed a similar algorithm in the case of windowed DFT
transforms. Inspired by the PhaseLift and MaxCut algorithms, but operating in the coefficients
space, the authors of \cite{W12} proposed a SemiDefinite Program called {\em PhaseCut}. They
show the algorithm yields the exact solution in the absence of noise under similar conditions
as PhaseLift.

The paper \cite{Bal12a} presents an iterative regularized least-square algorithm for
inverting the nonlinear map and compares its performance to a Cramer-Rao lower bound for this
problem in the real case. The paper also presents some new injectivity results which are
incorporated into this paper.

A different approach is proposed in \cite{ABFM12}. There the authors use a 4-term polarization identity together
with a family of spectral expander graphs to design a frame of bounded redundancy ($\frac{m}{n}\leq 236$) that
yields an exact reconstruction algorithm in the absence of noise.

The authors of \cite{EM12} study several robustness bounds to the phase recovery problem in the real case.
However their approach is different than ours in several respects. First they consider a probabilistic
setup of this problem, where data $x$ and frame vectors $f_j$'s are random vectors with probabilities
from a class of subgaussian  distributions. Additionally, their focus is on classes of $k$-sparse signals.
Their results show that, with high probability, recovery is possible from a number of measurements $m$ 
that has a similar asymptotic behavior with respect to $n$ and $k$ as in the case of linear measurements
(that is with the phase). In our paper we analyze stability bounds of reconstruction for a fixed frame
using deterministic analytic tools. After that we present asymptotic behavior of these bounds for random frames.

Finally, the authors of \cite{BCMN13a} analyze the phaseless reconstruction problem for both 
the real and complex case.  In the real case the authors obtain the exact upper Lipschitz 
constant for the nonlinear map $\nonlin$, namely $\sqrt{B}$ where $B$ is the upper frame bound. 
For the lower Lipschitz constant, they give an estimate between two computable singular 
eienvalues. Our results have overlaps with their results somewhat here. However, in this paper
we improve the improve the lower Lipschitz constant by giving its exact value. There are
some significant differences between this paper and \cite{BCMN13a}. In addition to 
studying of the Lipschitz property of the map $\nonlin$ we focus also on two related but
different settings. First we study the robustness of the reconstruction given a fixed error
allowance in measurements. Second we also consider the Lipschitz property of the map
$\nonlin^2$. The authors of \cite{BCMN13a} point out that the map $\nonlin^2$ is not 
bi-Lipschitz. However in our paper we show $\nonlin^2$ becomes bi-Lipschitz for a 
different metric on the domain. With this metric (the one induced by the nuclear norm
on the set of symmetric operators) the nonlinear map $\nonlin^2$ is bi-Lipschitz with 
constants indicated in
Theorem \ref{theo-4.5}. Furthermore the same conclusion holds true in the complex case, 
although this will be studied elsewhere.

The organization of the paper is as follows. Section \ref{sec2} formally defines the
problem and reviews existing inversion results in the real case.  Section \ref{sec3}
establishes information theoretic performance bounds, namely the Cramer-Rao lower bound.
Section \ref{sec4} contains robustness measures of any reconstruction algorithm. 
Section \ref{sec5} presents a stochastic analysis of these bounds, and is followed by references.

\section{Background}  \label{sec2}
\setcounter{equation}{0}

Let us denote by $H=\R^n$ the n-dimensional real Hilbert space $\R^n$ with scalar product
$\ip{}{}$. Let $\fc=\{f_1,\ldots,f_m\}$ be a spanning set of $m$ vectors in $H$. In finite
dimension (as it is the case here) such a set forms a {\em frame}. In the infinite dimensional
case, the concept of frame involves a stronger property than completeness (see for instance
\cite{Cass-artofframe}). We review additional terminology and properties which remain still
true in the infinite dimensional setting. The set $\fc$ is a frame if and only if there are two
positive constants $0<A\leq B<\infty$ (called frame bounds) so that
\begin{equation}  \label{2.1}
     A\norm{x}^2 \leq \sum_{k=1}^m |\ip{x}{f_k}|^2 \leq B\norm{x}^2.
\end{equation}
When we can choose $A=B$ the frame is said {\em tight}. For $A=B=1$ the frame is called {\em
Parseval}. The {\em frame matrix} corresponding to $\fc$ is defined as
$F=[f_1,f_2, \dots, f_m]$ with the vectors $f_j\in\fc$ as its columns.
We shall frequently identify $\fc$ with its corresponding frame matrix $F$. The largest
$A$ and smallest $B$ in \eqref{2.1} are called the {\em lower frame bound} and
{\em upper frame bound} of $\fc$, and they are given by
\begin{equation}  \label{2.2}
    A = \lambda_{\mathrm{max}}(FF^*) =\sigma_1^2(F), \mhsp
    B = \lambda_{\mathrm{min}}(FF^*) =\sigma_n^2(F)
\end{equation}
where $\lambda_{\mathrm{max}}, \lambda_{\mathrm{min}}$ denote the largest and smallest
eigenvalues respectively, while $\sigma_1, \sigma_n$ denote the first and $n$-th singular
values respectively.
A set of vectors $\fc$ of the $n$-dimensional Hilbert space $H$ is said to be
{\em full spark} if any subset of $n$ vectors is linearly independent.

For a vector $x\in H$, the collection of coefficients $\{\ip{x}{f_j}:~1\leq j\leq m\}$
represents the analysis map of vector $x$ given by the frame $\fc$, and from which
$x$ can be completely reconstructed. In the phaseless reconstruction problem,
we ask the following question: Can $x$ be reconstructed from
$\{|\ip{x}{f_j}|:~1\leq j\leq m\}$? Consider the following equivalence relation
$\sim$ on $H$: $x\sim y$ if and only if $y=cx$ for some unimodular
constant $c$, $|c|=1$. Since we focus on the real vector space $H=\R^n$, we have
$x \sim y$ if and only if $x=\pm y$. Clearly the phaseless reconstruction problem cannot
distinguish $x$ and $y$ if $x\sim y$, so we will be looking at reconstruction
on $\hat{H}:=H/\sim = \R^n/\sim$ whose elements are given by equivalent classes
$\hat{x}=\{x,-x\}$ for $x\in\R^n$. The analogous analysis map for phaseless reconstruction
is the following nonlinear map
\begin{equation} \label{eq:nonlin}
   \nonlin :\hat{H}\rightarrow \R_+^m, \mhsp
   \nonlin(\hat{x})=[|\ip{x}{f_1}|, |\ip{x}{f_2}|, \dots, |\ip{x}{f_m}|]^T.
\end{equation}
Note that $\nonlin$ can also be viewed as a map from $\R^n$ to $\R_+^m$. Throughout the paper
we will not make an explicit distinction unless such a distinction is necessary.

Thus the phaseless reconstruction problems aims to reconstruct $\hat x\in \hat H$ from
the map $\nonlin(x)$. We say a frame $\fc$ is {\em phase retrievable} if one can
reconstruct $\hat x\in \hat H$ for all $\hat x$, or in other words,  
$\nonlin$ is injective on $\hat H$.
The main objective of this paper is to analyze robustness and stability of the
inversion map, and to give
performance bounds of any reconstruction algorithm.

Before proceeding further we first review existing results on injectivity of
the nonlinear map $\nonlin$. In general a subset
$Z$ of a topological space 
is said {\em generic} if its open interior is dense. However in the following statements, the term
{\em generic} refers to Zarisky topology: a set $Z\subset \K^{n\times m}= \K^{n}\times\cdots\times \K^n$ is said {\em generic}
if $Z$ is dense in $\K^{n\times m}$ and its complement is a finite union of zero sets of polynomials in $nm$ variables
with coefficients in the field $\K$ (here $\K=\R$).

\begin{thm} \label{theo-2.1}
Let $\fc$ be a frame in $H=\R^n$ with $m$ elements. Then the following hold true:
\begin{enumerate}
\item The frame $\fc$ is phase retrievable in $\hat H$
   if and only if for any disjoint
   partition of the frame set $\fc=\fc_1\cup\fc_2$, either $\fc_1$ spans $\R^n$ or
   $\fc_2$ spans $\R^n$.

\item If $\fc$ is phase retrievable in $\hat H$ then $m\geq 2n-1$.
Furthermore, for a generic $\fc$ with $m \geq 2n-1$ the map $\nonlin$
is phase retrievable in $\hat H$.

\item Let $m=2n-1$. Then  $\fc$ is phase retrievable in $\hat H$
if an only if $\fc$ is full spark.

\item Let 
\begin{equation}
\label{eq:a0}
a_0 := \min_{\|x\|=\|y\|=1}\sum_{j=1}^m |\ip{x}{f_j}|^2|\ip{y}{f_j}|^2  \geq 0.
\end{equation}
Then
\begin{equation}
\label{eq:q4}
   \sum_{k=1}^m |\ip{x}{f_k}|^2 |\ip{y}{f_k}|^2 \geq a_0 \norm{x}^2 \norm{y}^2.
\end{equation}
Then $\fc$ is phase retrievable in $\hat H$ if and only if $a_0>0$.

\item For any $x\in\R^n$ define the matrix $R(x)$ by
\begin{equation} \label{eq:q2}
    R(x):= \sum_{j=1}^m |\ip{x}{f_j}|^2 f_jf_j^*.
\end{equation}
Then $R(x) \geq a_0\|x\| I$ where $I$ is the identity matrix and $a_0$ is given by (\ref{eq:a0}). In other words,
$\lambda_{\mathrm{min}}(R(x)) \geq a_0\|x\|^2$. Furthermore $a_0=\min_{\|x\|=1}\lambda_{\mathrm{min}}(R(x))$.

\end{enumerate}
\end{thm}
\Proof The results (1)-(3) are in \cite{BCE06}, and (4)-(5) are in  \cite{Bal12a}.
\eproof

\section{Information Theoretic Performance Bounds\label{sec3}}
\setcounter{equation}{0}

In this section we derive expressions for the {\em Fisher Information Matrix}
and obtain performance
bounds for reconstruction algorithms in the noisy case.

Consider the following noisy measurement process:
\begin{equation}
\label{eq:noise1}
y_k = |\ip{x}{f_k}|^2 + \nu_k ~~,~~\nu_k\sim \Ncal(0,\sigma^2)~,~1\leq k\leq m
\end{equation}
where the noise model is AWGN (additive white Gaussian noise): each random variable $\nu_k$ is independent and normally distributed with
zero mean and $\sigma^2$ variance.

Consider the noiseless case first (that is $\nu_k=0$).
Obviously one cannot obtain the exact vector $x\in H$ due to the global sign ambiguity.
Instead the best outcome is to identify (that is, to estimate) the class $\hat{x}=\{x,-x\}$ from $\nonlin(x)$. As such, we fix a disjoint partition of the punctured
Hilbert space $H$, $\R^n\setminus\{0\}=\Omega_1\cup\Omega_2$, such that $\Omega_2=-\Omega_1$. We make the choice that the vector $x$ belongs to
$\Omega_1$. Hence any estimator of $x$ is a map
$\recon:\R^m\longrightarrow \Omega_1\cup\{0\}$. Denote by $\intOmega$ its interior
as a subset of $\R^n$. A typical such decomposition is
$$
   \Omega_1=\bigcup_{k=1}^{n}\Bigl\{x\in\R^n:~x_k \geq 0,~x_j =0 ~
         \mbox{$x_j=0$ for $j<k$}\Bigr\}.
$$
Note its interior is given by $\intOmega=\{x\in\R^n~,~x_1>0\}$.

Under these assumptions we compute the Fisher Information matrix (see \cite{Kay2010}).
This is given by
\begin{equation} \label{eq:Fish}
 (\I(x))_{k,j} = \E \left[ (\nabla \log\,L(x)) (\nabla \log\,L(x))^T \right]
\end{equation}
where the likelihood function $L(x)$ is given by
\begin{equation}
\label{eq:like}
L(x)=p(y|x) = \frac{1}{(2\pi)^{m/2}\sigma^m}
    \exp\bigl(-\frac{1}{2\sigma^2}\norm{y-\nonlin(x)}^2 \big).
\end{equation}
After some algebra (see \cite{Bal12a}) we obtain
\begin{equation} \label{eq:Fisher}
\I(x) = \frac{4}{\sigma^2} R(x)~~,~~R(x)=\sum_{j=1}^m |\ip{x}{f_j}|^2 f_jf_j^T.
\end{equation}
Note the matrix $R(x)$ is exactly the same as the matrix introduced in (\ref{eq:q2}). Thus we obtain the following results:

\begin{thm} \label{th3.1}
The frame $\fc$ is phase retrievable if and only if the
Fisher information matrix $\I(x)$ is invertible for any $x\neq 0$.
Furthermore, when $\fc$ is phase retrievable there is a positive constant $a_0>0$ so that
\begin{equation}
\I(x) \geq \frac{4 a_0}{\sigma^2} \norm{x}^2 I
\end{equation}
where $I$ is the $n\times n$ identity.
\end{thm}

This allows to state the following performance bound result (see \cite{Kay2010}
for details on the Cramer-Rao lower bound).

\begin{thm}\label{th3.2}
   Assume $x\in\intOmega$.  Let $\recon:\R^m\rightarrow\Omega_1$ be any unbiased estimator
   for $x$. Then its covariance matrix is bounded below by the Cramer-Rao lower bound:
  \begin{equation} \label{eq:CRLB}
     \mathrm{Cov}[\recon(y)] \geq (\I(x))^{-1} =\frac{\sigma^2}{4} (R(x))^{-1}.
  \end{equation}
   Furthermore, any efficient estimator (that is, any unbiased estimator $\recon$ that
   achieves the Cramer-Rao Lower Bound (\ref{eq:CRLB})) has the covariance matrix bounded from above by
   \begin{equation}
        \mathrm{Cov}[\recon(y)] \leq \frac{\sigma^2}{4a_0 \norm{x}^2 } I
   \end{equation}
   and Mean-Square error bounded above by
   \begin{equation} \label{eq:MSE}
      \mathrm{MSE}(\recon) = \E \left[ \norm{\recon(y)-x}^2 \right] \leq \frac{n\sigma^2}{4a_0 \norm{x}^2}.
   \end{equation}
\end{thm}

\section{Robustness Measures for Reconstruction}
\label{sec4}
\setcounter{equation}{0}

In this section we analyze the robustness of deterministic phaseless reconstruction.
Additionally we connect the constant $a_0$ introduced earlier in 
Theorem \ref{theo-2.1} to quantities directly
computable from the frame $\fc$.

A natural approach is to analyze the stability in the worst case scenario, for which we consider the
following measures. Denote $d(x,y) := \min (\norm{x-y},\norm{x+y})$. For any $x\in \R^n$ define
\begin{equation} \label{4.1}
    \qq_\eps (x) =  \max_{\{y:\norm{\nonlin(x)-\nonlin(y)} \leq \eps\}} \frac{d(x,y)}{\eps}.
\end{equation}
The size of $\qq_\eps(x)$ measures the worst case
stability of the reconstruction for the vector $x$, under the assumption that the
total noise level is controlled by $\eps$. We also
study the global stability by analyzing the measures
\begin{equation}  \label{4.2}
q_\eps := \max_{\norm{x}=1} \qq_\eps (x), \shsp
q_0 := \limsup_{\eps\rightarrow 0}\, q_\eps, \shsp
q_{\infty} := \sup_{\eps>0}\, q_\eps.
\end{equation}
Here $\|.\|$ denotes usual Euclidian norm. Note that $\qq_\eps(x)$ has the scaling
property $\qq_\eps(x) = \qq_{|c|\eps}(cx)$ for any real $c\neq 0$. Thus it is natural
to focus on unit vectors $x$.

We introduce now some quantities that play key roles in the estimation of these robustness measures.
For the frame $\fc$ let $F=[f_1, f_2, \cdots, f_m]$ be its frame matrix.
Denote by $\fc[S]=\{f_k~,~k\in S\}$ the subset of $\fc$
indexed by a subset $S \subseteq \{1,2,\ldots,m\}$, and
by $F_S$ the frame matrix corresponding to $\fc[S]$ (which is the matrix with vectors in $\fc[S]$
as its columns). Set
\begin{equation}  \label{4.3}
     A[S] := \sigma_n^2(F_S) = \lambda_{\mathrm{min}}(F_SF_S^*),
\end{equation}
where as usual $\sigma_n$ and $\lambda_{\mathrm{min}}$ denote the $n$-th singular value and the minimal eigenvalue, respectively.
Note that $A[S]$ is in fact the lower frame bound of $\fc[S]$.

Let $\Sc$ denote the collection of subsets $S$ of $\{1,2,\ldots,m\}$ so that
$\dim\,(\mathrm{span}(\fc[S^c]))<n$, where $S^c=\{1,2,\ldots,m\}\setminus S$ is the complement of $S$.
In other words, $\rank(F_{S^c}) <n$.
Denote by $\Delta$ and $\omega$ the following expressions:
\begin{eqnarray}
\Delta & = & \min_S\,\sqrt{A[S]+A[S^c]} \label{eq:Delta} \\
\omega & = & \min_{S\in\Sc}\,\sigma_n(F_S). \label{eq:omega}
\end{eqnarray}
All of them depend of course on $\fc$. However since we fix $\fc$ throughout the paper, we shall
without confusion
not explicitly reference $\fc$ in the notation for simplicity as there will not be any confusion.
Clearly
\begin{equation} \label{4.4}
   \Delta \leq \omega.
\end{equation}

\begin{proposition}  \label{prop-4.1}
Let $\eps>0$. Then the stability measurement function $\qq_\eps(x)$ is given by
\begin{equation} \label{4.5}
    \qq_\eps(x) = \frac{1}{\eps} \max_{w_1,w_2} \min\bigl\{\|w_1\|, \|w_2\|\bigr\}
\end{equation}
under the constraints $\frac{1}{2} (w_1+w_2)=x$ and
\begin{equation} \label{4.6}
    \sum_{j=1}^m \min\,\bigl(|\inner{f_j,w_1}|^2, |\inner{f_j,w_2}|^2\bigr)=
    \bigl\|F_S^* w_1\bigr\|^2 + \bigl\|F_{S^c}^* w_2\bigr\|^2 \leq \eps^2,
\end{equation}
where $S := S(w_1,w_2)=\{j:~|\inner{f_j,w_1}| \leq |\inner{f_j,w_2}|\}$.
\end{proposition}
\Proof For any $x, y\in\R^n$ let $w_1=x+y$ and $w_2=x-y$. Then $x=\frac{1}{2}(w_1+w_2)$
and $y=\frac{1}{2}(w_1+w_2)$. It is easy to check that for
$S =\{j:~|\inner{f_j,w_1}| \leq |\inner{f_j,w_2}|\}$ we have
$$
    |\inner{f_j,x}| - |\inner{f_j,y}| = \left\{\begin{array}{cl}
         \pm \inner{f_j,w_1} \mhsp& j\in S, \\
         \pm \inner{f_j,w_2} \mhsp& j\in S^c. \end{array}\right.
$$
In other words,
\begin{equation} \label{4.7}
   |\inner{f_j,x}| - |\inner{f_j,y}| = \min (|\inner{f_j,w_1}|, |\inner{f_j,w_2}|).
\end{equation}
Let $F$ be the frame matrix of $\fc$. We thus have
$$
    \Bigl\| \nonlin(x)-\nonlin(y) \Bigr\|^2 = \sum_{j\in S} |\inner{f_j,w_1}|^2 +
          \sum_{j\in S^c} |\inner{f_j,w_2}|^2
          = \bigl\|F_S^* w_1\bigr\|^2 + \bigl\|F_{S^c}^* w_2\bigr\|^2.
$$
Note that $d(x,y) = \min (\|w_1\|, \|w_2\|)$. The proposition now follows.
\eproof

\vspace{2mm}

The above proposition allows us to establish the following stability result for the worst case scenario.

\begin{thm}  \label{theo-4.2}
   Assume that the frame $\fc$ is phase retrievable. Let $A>0$ be the lower frame bound for
   the frame $\fc$ and let
   $\tau:=\min \{\sigma_n(F_S):~S\subseteq \{1,\dots,m\}, ~\rank(F_S) = n\}$.
   \begin{itemize}
   \item[\rm (A)]~~For any $\eps>0$ we have
   \begin{equation} \label{4.8}
    \min\Bigl\{\frac{1}{\eps}, \frac{1}{\omega}\Bigr\} \leq q_\eps \leq \frac{1}{\Delta}.
   \end{equation}
   \vspace{2mm}
   \item[\rm (B)]~~If $\eps < \tau$ then $q_\eps = \frac{1}{\omega}$.
   Consequently $q_0 = \frac{1}{\omega}$. \smallskip
   \item[\rm (C)]~~For any nonzero $x\in\R^n$ and any $0<\eps<\delta_x$ we have
   \begin{equation}  \label{4.9}
       \qq_\eps(x) = \frac{1}{\sqrt A},
   \end{equation}
   where
   $$
      \delta_x := \frac{2\tau}{\max(\|f_j\|) +\tau} \min\,
              \Bigl\{|\inner{f_j,x}|:~\inner{f_j,x} \neq 0\Bigr\}.
   $$
   \item[\rm (D)]~~The upper bound $q_\infty$ equals the reciprocal of $\Delta$:
   \begin{equation} \label{eq:4.8bis}
    q_\infty = \frac{1}{\Delta}.
   \end{equation}
   \end{itemize}
\end{thm}
\Proof  To prove (A) we first establish the upper bound in \eqref{4.8}. Let $x\in\R^n$.
By Proposition \ref{prop-4.1} we have
$$
    \qq_\eps(x) = \frac{1}{\eps} \max_{w_1,w_2} \min\bigl\{\|w_1\|, \|w_2\|\bigr\}
$$
under the constraints $\frac{1}{2} (w_1+w_2)=x$ and
$$
    \bigl\|F_S^* w_1\bigr\|^2 + \bigl\|F_{S^c}^* w_2\bigr\|^2 \leq \eps^2
$$
for some $S$. Now assume without loss of generality that $\|w_1\| \leq \|w_2\|$.
Then
\begin{eqnarray*}
    \frac{\eps^2}{\|w_1\|^2} &\geq&
    \frac{\bigl\|F_S^* w_1\bigr\|^2 + \bigl\|F_{S^c}^* w_2\bigr\|^2}{\|w_1\|^2}\\
    &\geq & \sigma_n^2(F_S) + \sigma_n^2(F_{S^c}) \frac{\|w_2\|^2}{\|w_1\|^2} \\
    &\geq& \Delta.
\end{eqnarray*}
It follows that
$$
    \frac{1}{\eps}\min\bigl\{\|w_1\|, \|w_2\|\bigr\} \leq \frac{1}{\Delta}.
$$
Thus $\qq_\eps(x) \leq \frac{1}{\Delta}$.

To establish the lower bound in \eqref{4.8} we construct for any $\eps>0$
an $x\in\R^n$ and vectors
$w_1, w_2$ satisfying the imposed constraints. Let $S$ be a subset of $\{1,2, \dots, m\}$ such that
$\rank(F_{S^c}) <n$ and $\sigma_n(F_S)=\omega$. Choose $v_1,v_2\in \R^n$ with the property
$\norm{v_1}=\norm{v_2}=1$ and
$$
   \| F^*_Sv_1\| = \omega, \mhsp F^*_{S^c} v_2=0.
$$
Set
$$
   t = \min\Bigl\{\frac{\eps}{\omega}, 1\Bigr\}, \shsp \mbox{and} \shsp w_1=tv_1.
$$
Hence $\norm{w_1}=t\leq 1$. Now we select an $s\in\R$ so that $\norm{w_1+sv_2}=2$. This is always
possible since $s\mapsto \norm{w_1+sv_2}$ is continuous and
$\norm{w_1+0v_2}=t\leq 1\leq 2 \leq \norm{w_1+3v_2}$.
Set $w_2=sv_2$ so $\norm{w_1+w_2}=2$.
We have
$$
     |s|=\norm{sv_2}\geq \norm{w_1+sv_2}-\norm{w_1} = 2-t\geq 1.
$$
Thus $\norm{w_2}\geq \norm{w_1}$. Now let
$$
    x=\frac{1}{2}(w_1+w_2) \mhsp \mbox{and} \mhsp y=\frac{1}{2}(w_1-w_2).
$$
We have then
\begin{eqnarray*}
   \norm{\nonlin(x)-\nonlin(y)}^2 &=&
   \sum_{j=1}^m \min\,(|\inner{f_j,w_1}|^2,|\inner{f_j,w_2}|^2) \\
   &\leq& \sum_{j\in S}|\inner{f_j,w_1}|^2 + \sum_{j\in S^c}|\inner{f_j,w_2}|^2 \\
   &=& t^2 \omega^2 ~\leq ~\eps^2.
\end{eqnarray*}
Furthermore
$$
   d(x,y) = \min(\norm{w_1},\norm{w_2}) = \norm{w_1} = t.
$$
Hence for this $x$ we have
$$
    \qq_\eps(x) \geq \frac{d(x,y)}{\eps} = \min \Bigl\{\frac{1}{\eps}, \frac{1}{\omega}\Bigr\}.
$$
It follows that $q_\eps \geq \min\,\{\frac{1}{\eps}, \frac{1}{\omega}\}$. Now by
taking $\eps>0$ sufficiently small we have $q_\eps \geq \frac{1}{\omega}$. 


We now prove (B).
Assume that $\eps \leq \min \{\sigma_n(F_S):~\rank(F_S) = n\}$. Then clearly we have
$\eps \leq \omega$. Thus by \eqref{4.8} we have $q_\eps \geq \frac{1}{\omega}$.
Again for each $x\in\R^n$ with $\|x\|=1$ we consider $w_1,w_2$ for the estimation of  $q_\eps(x)$.
The constraint $\|w_1+w_2\|=2$ implies either $\|w_1\|\geq 1$ or $\|w_2\|\geq 1$. Without
loss of generality we assume that $\|w_1\|\geq 1$.  For the constraint
$\bigl\|F_S^* w_1\bigr\|^2 + \bigl\|F_{S^c}^* w_2\bigr\|^2 \leq \eps^2$ for some $S$, assume that
$\rank(F_S) =n$ then we have
$$
  \|F^*_Sw_1\| \geq \sigma_n(F_S) \|w_1\| \geq \min \{\sigma_n(F_S):~\rank(F_S) = n\}> \eps.
$$
This is a contradiction. So $\rank(F_S) <n$ and hence
$$
   \eps^2 \geq \bigl\|F_S^* w_1\bigr\|^2 + \bigl\|F_{S^c}^* w_2\bigr\|^2
    \geq \bigl\|F_{S^c}^* w_2\bigr\|^2\geq \omega^2\|w_2\|^2.
$$
Thus $\|w_2\| \leq \frac{\eps}{\omega}$. Proposition \ref{prop-4.1} now yields
$q_\eps = \frac{1}{\omega}$, proving part (B).

Now we prove (C). We go back to the formulation in Proposition \ref{prop-4.1}.
$$
    \qq_\eps(x) = \frac{1}{\eps} \max_{w_1,w_2} \min\bigl\{\|w_1\|, \|w_2\|\bigr\}
$$
under the constraints $\frac{1}{2} (w_1+w_2)=x$ and
$$
    \bigl\|F_S^* w_1\bigr\|^2 + \bigl\|F_{S^c}^* w_2\bigr\|^2 \leq \eps^2
$$
where $S := S(w_1,w_2)=\{j:~|\inner{f_j,w_1}| \leq |\inner{f_j,w_2}|\}$. Since
$\nonlin$ is injective, either $\rank(F_S)=n$ or $\rank(F_{S^c})=n$ by Theorem \ref{theo-2.1}
(1). Without loss of generality we assume $\rank(F_S)=n$. Thus
$\eps\geq \|F_S^* w_1\bigr\| \geq \tau \|w_1\|$. So
$\|w_1\| \leq \eps/\tau$. We show that for any $k\in S^c$ we must
have $\inner{f_k,x}=0$. Assume otherwise and write $w_2=2x-w_1$,
$L_x:= \min\,\{|\inner{f_j,x}|:~\inner{f_j,x} \neq 0\}$. Then
$$
  |\inner{f_k, w_2}| \geq 2|\inner{f_k,x}|-|\inner{f_k,w_1}|
  \geq 2 L_x - \max\,(\|f_j\|) \|w_1\| \geq 2 L_x -\max\,(\|f_j\|)\frac{\eps}{\tau}>\eps.
$$
This is a contradiction. Thus for $k\in S^c$ we have $\inner{f_k,x}=0$ and
$$
   |\inner{f_j,w_2}| = |\inner{f_j,2x-w_1}| =|\inner{f_j,w_1}|.
$$
It follows that
$$
    \bigl\|F_S^* w_1\bigr\|^2 + \bigl\|F_{S^c}^* w_2\bigr\|^2 =\|F^* w_1\|^2\leq \eps^2.
$$
Thus $\|w_1\| \leq \eps/\sqrt A$ and hence $\qq_\eps(x) \leq \frac{1}{\sqrt A}$.
Now we show the bound can be achieved. Let $w_1$ satisfy $\|F^*w_1\| = \sqrt A \|w_1\|
=\eps$. Such a $w_1$ always exists. Then clearly $w_1$ and $w_2 = 2x-w_1$ satisfy
the required constraints, and it is easy to check that
$\min\,(\|w_1\|, \|w_2\|) = \|w_1\| = \eps/\sqrt A$.

Finally we prove (D). By the result at part (A),  $q_\infty \leq \frac{1}{\Delta}$.
It is therefore sufficient to shoow that $Q_\eps(x)\geq \frac{1}{\Delta}$ for some $x$ and $\eps$.
Let $S_0$ be the subset that achieves the minimum in (\ref{eq:Delta}). Let $u,v\in H$ be unit eigenvectors
corresponding to the lowest eigenvalues of $F_{S_0}F_{S_0}^*$ and $F_{S_0^c}F_{S_0^c}^*$
respectively. Thus
\[ \norm{F_{S_0}^*u}^2 = A[S_0]~~,~~\norm{F_{S_0^c}^*v}^2 = A[S_0^c] \]
Let $x=(u+v)/2$ and $\eps=\Delta$, and set $w_1=u$, $w_2=v$. Then by Proposition \ref{prop-4.1}
\[ Q_\eps(x) \geq \frac{\min(\norm{w_1},\norm{w_2})}{\eps} = \frac{1}{\Delta} \]
since
\[ \sum_{j=1}^m \min(|\ip{f_j}{w_1}|^2,|\ip{f_j}{w_2}|^2) \leq \norm{F_{S_0}^*w_1}^2
+ \norm{F_{S_0^c}^* w_2}^2 = \eps^2 \]
This concludes the proof.
\eproof

\vspace{3mm}
\noindent
{\bf Remark.}~~It may seem strange that $\qq_\eps(x) = \frac{1}{\sqrt A}$ for all $x\neq 0$
and sufficiently small $\eps$ while $q_0 = \frac{1}{\omega}$, where $\omega$ is typically
much smaller than $\sqrt A$. The reason is that for $\qq_\eps(x) = \frac{1}{\sqrt A}$ to
hold, $\eps$ depends on $x$. Thus we cannot exchange the order of $\limsup_{\eps\ra0}$
and $\max_{\|x\|=1}$.

\bigskip

Related to the study of stability of phaseless reconstruction is the study
of Lipschitz property of the map $\nonlin$ on $\hat H :=\R^n/\sim$.
We analyze the bi-Lipschitz bounds of both $\nonlin$ and $\nonlin^2$, which is
simply the map $\nonlin$ with all entries squared, i.e.
$$
    \nonlin^2(x) := [|\inner{f_j,x}|, \dots, |\inner{f_m,x}|]^T.
$$
We shall consider two distance functions on $\hat H =\R^n/\sim$: the standard
distance $d(x,y) :=\min (\|x-y\|, \|x+y\|)$ and the distance $d_1(x,y):=
\|xx^*-yy^*\|_1$ where $\|X\|_1$ denotes the {\em nuclear norm} of $X$, which is the
sum of all singular values of $X$.
Specifically we are interested in examining the local and global behavior of the ratios
\begin{equation} \label{4.10}
    U(x,y) : = \frac{\|\nonlin(x) - \nonlin(y)\|}{d(x,y)}, \mhsp
    V(x,y) : = \frac{\|\nonlin^2(x) - \nonlin^2(y)\|}{d_1(x,y)}.
\end{equation}

\vspace{2mm}

We first investigate the bounds for $U(x,y)$. For this the upper bound is relatively straightforward.
Let $w_1=x-y$ and $w_2=x+y$. We have already shown in the proof of Theorem \ref{theo-4.2} using
\eqref{4.7} that
\begin{eqnarray*}
   \norm{\nonlin(x)-\nonlin(y)}^2 &=&
   \sum_{j=1}^m \min\,(|\inner{f_j,w_1}|^2,|\inner{f_j,w_2}|^2) \\
   &\leq& \min\,\Bigl\{\sum_{j=1}^m |\inner{f_j,w_1}|^2, \sum_{j=1}^m |\inner{f_j,w_2}|^2\Bigr\} \\
   &\leq& B \min\, \Bigl\{\|w_1\|^2, \|w_2\|^2\Bigr\} ~=~ B d^2(x,y),
\end{eqnarray*}
where $B$ is the upper frame bound of the frame $\fc$. Thus
$U(x,y)$ has an upper bound $U(x,y) \leq \sqrt B$. Furthermore, the bound is sharp.
To see this, pick a unit vector $x\in\R^n$ such that
$\sum_{j=1}^m |\inner{f_j,w_1}|^2=B$ and set $y=2x$. Then $U(x,y) = \sqrt B$.

To study the lower bound $U(x,y)$ we now consider the following quantities:
\begin{eqnarray*}
   \rho_\eps (x)  &:=&  \inf_{\{y: d(x,y)\leq \eps\}} U(x,y), \\
   \rho(x)  &:=&  \liminf_{\{y: d(x,y)\ra 0\}} U(x,y) = \liminf_{\eps\ra 0} \rho_\eps(x), \\
   \rho_0 &:=&  \inf_x \rho(x) , \\
   \rho_{\infty}  &:=&  \inf_{d(x,y)>0} U(x,y).
\end{eqnarray*}
We apply the equality
$$
    U^2(x,y) = \frac{\sum_{j=1}^m \min\,(|\inner{f_j,w_1}|^2,|\inner{f_j,w_2}|^2)}
               {\min\,(\norm{w_1}^2,\norm{w_2}^2)}
$$
where again $w_1=x-y$ and $w_2=x+y$.
Now fix $x$ and let $d(x,y) < \eps$. Without loss of generality we may assume
$\|y - x\|<\eps$. Thus $\|w_1\|<\eps$ and $\|w_2- 2x\|=\|w_1\|<\eps$.
Let $S = \{j,~\inner{f_j,x}\neq 0\}$ and set
\begin{equation} \label{4.11}
 \eps_0(x) :=\frac{ \min_{k\in S}|\inner{f_k,x}|}{\max_{k\in S} \norm{f_k}}.
\end{equation}
Note for any $w_1$ with $\norm{w_1}<\eps_0$ and $k\in S$ we have
$$
   |\inner{f_k,w_2}| =|2\inner{f_k,x}-\inner{f_k,w_1}|
    \geq 2|\inner{f_k, x}|-|\inner{f_k, w_1}|
     \geq 2\eps_0(x)\norm{f_k}-\norm{w_1}\norm{f_k}\geq |\inner{f_k,w_1}|,
$$
whereas for $k\in S^c$ we have
$$
    |\inner{f_k, w_2}|=|\inner{f_k, w_1}|.
$$
Hence $\min\,(|\inner{f_j,w_1}|^2,|\inner{f_j,w_2}|^2)=|\inner{f_j,w_1}|^2$
for all $j$ whenever $\eps < \eps_0(x)$. It follows that
$$
    U^2(x,y) = \frac{\sum_{j=1}^m|\inner{f_j,w_1}|^2}{\norm{w_1}^2}
    =\sum_{j=1}^m\bigl|\ip{\frac{w_1}{\norm{w_1}}}{f_j}\bigr|^2.
$$
Thus $U^2(x,y) \geq A$ where $A$ is the lower frame bound for the frame $\fc$. Furthermore
this lower bound is achieved whenever $w_1=x-y$ is an eigenvector corresponding to
the smallest eigenvalue of $FF^*$. This implies that
$$
    \rho_\eps(x) = \sqrt A
$$
whenever $\eps<\eps_0(x)$. Consequently $\rho(x) = \sqrt A$.
We have the following theorem:

\begin{thm} \label{theo-4.3}
    Assume that the frame $\fc$ is phase retrievable. Let $A$, $B$ be the lower and upper
frame bounds for the frame $\fc$, respectively and for each $x \in\R^n$, let
$\eps_0(x)$ be given in \eqref{4.11} . Then
    \begin{itemize}
    \item[\rm (1)]~~$U(x,y) \leq \sqrt B$ for any $x,y\in\R^n$ with $d(x,y)>0$.
    \item[\rm (2)]~~Assume that $\eps < \eps_0(x)$. Then $\rho_\eps(x) = \sqrt A$.
       Consequently $\rho(x) = \rho_0 = \sqrt A$.
    \item[\rm (3)]~~$\Delta = \rho_{\infty} \leq \omega\leq \rho_0 =\rho(x) = \sqrt{A}$.
    \item[\rm (4)]~~ The map $\nonlin$ is bi-Lipschitz with (optimal) upper Lipschitz bound $\sqrt{B}$ and lower Lipschitz bound 
$\rho_{\infty}$:
\[ \rho_{\infty}  d(x,y) \leq \norm{\nonlin(x)-\nonlin(y)} \leq \sqrt{B} d(x,y) ~~,~~\forall x,y\in \hat{H} \]
    \end{itemize}
\end{thm}
\Proof   We have already proved (1) and (2) of the theorem in the discussion. It remains only
to prove (3) since (4) is just a restatement of (1) and (3). Note that
$$
 \rho^2_{\infty}  =  \inf_{d(x,y)>0} U^2(x,y)
         =\inf_{w_1, w_2\neq 0} \frac{\sum_{j=1}^m \min\,(|\ip{f_j}{w_1}|^2,|\ip{f_j}{w_2}|^2)}
                               {\min\,(\norm{w_1}^2,\norm{w_2}^2)} .
$$
For any $w_1, w_2$, assume without loss of generality that $0<\|w_1\| \leq \|w_2\|$. Let
$S=\{j:~|\ip{f_j}{w_1}|\leq |\ip{f_j}{w_2}|\}$. Set $v_1=w_1/\|w_1\|$, $v_2=w_2/\|w_2\|$ and
$t=\|w_2\|/\|w_1\|\geq 1$. Then
\begin{eqnarray*}
 \frac{\sum_{j=1}^m \min\,(|\ip{f_j}{w_1}|^2,|\ip{f_j}{w_2}|^2)} {\min\,(\norm{w_1}^2,\norm{w_2}^2)}
  &=& \sum_{j\in S} |\ip{f_j}{v_1}|^2 +t^2 \sum_{j\in S^c} |\ip{f_j}{v_2}|^2\\
&\geq & \sum_{j\in S} |\ip{f_j}{v_1}|^2 + \sum_{j\in S^c} |\ip{f_j}{v_2}|^2\\
&\geq & \Delta^2.
\end{eqnarray*}
Hence $\rho_\infty \geq \Delta$.

\ignore{
Observe that the above allows us to rewrite $\rho^2_\infty$ as
\begin{equation}  \label{4.12}
 \rho^2_{\infty} =\inf_{\|v_1\|=\|v_2\|=1\leq t}
    \sum_{j=1}^m \min\,(|\ip{f_j}{v_1}|^2, t^2|\ip{f_j}{v_2}|^2).
\end{equation}
Again we consider the set $\Sc:=\{S\subseteq \{1, \dots,m\}:~\rank(F_{S^c})<n\}$. For any
$S\in\Sc$ set
$$
    J_S:=\Bigl\{(v_1,v_2,t)\in\R^n\times\R^n\times\R:~\|v_1\|=\|v_2\|=1\leq t,F_{S^c}^*v_2=0\Bigr\}.
$$
Then taking infimum over a smaller set in \eqref{4.12} we have
\begin{eqnarray*}
\rho^2_{\infty} &\leq& \inf_{S\in\Sc}\,\inf_{(v_1,v_2,t)\in J_S}
   \sum_{j=1}^m \min\,(|\ip{f_j}{v_1}|^2, t^2|\ip{f_j}{v_2}|^2) \\
&= & \inf_{S\in \Sc}\,\inf_{\|v_1\|=\|v_2\|\leq t} \sum_{j\in S}
    \min\,(|\ip{f_j}{v_1}|^2, t^2|\ip{f_j}{v_2}|^2) \\
&\leq& \inf_{S\in\Sc}\,\inf_{\|v_1\|=1} \sum_{j\in S} |\ip{f_j}{v_1}|^2 \\
&=& \omega^2.
\end{eqnarray*}
}

Let $S$ and $u,v\in H$ be normalized (eigen) vectors that achieve the bound $\Delta$, that is:
$$
  \norm{u}=\norm{v}=1,\mhsp 
  \sum_{k\in S}|\ip{u}{f_k}|^2 + \sum_{k\in S^c}|\ip{v}{f_k}|^2 = \Delta.
$$
Set $x=u+v$ and $y=u-v$. Then, following \cite{BCMN13a}
\begin{eqnarray*}
 \norm{\nonlin(x)-\nonlin(y)}^2 & = & \sum_{k\in S} \left|\, |\ip{u+v}{f_k}|-|\ip{u-v}{f_k}|\, \right|^2 +
\sum_{k\in S^c} \left|\, |\ip{u+v}{f_k}|-|\ip{u-v}{f_k}| \,\right|^2  \\
& \leq & 4\left( \sum_{k\in S}|\ip{u}{f_k}|^2
+\sum_{k\in S^c}|\ip{v}{f_k}|^2 \right) = 4\Delta. 
\end{eqnarray*}
On the other hand
$$
   d(x,y) = min(\norm{x-y},\norm{x+y}) = 2.
$$
Thus we obtain
$$
   \frac{\norm{\nonlin(x)-\nonlin(y)}}{d(x,y)} \leq \Delta.
$$

The theorem is now proved.
\eproof

\vspace{3mm}
\noindent
{\bf Remark.}~~The two quantities, $\rho_\infty$ and $q_\infty$ satisfy $\rho_{\infty}= \frac{1}{q_\infty}$.
However there are subtle differences between $q_\eps(x)$ and $\rho_\eps(x)$ so that the simple relationship
$\rho_\eps(x) = 1/q_\eps(x)$ does not usually hold. 

\vspace{3mm}

\noindent
{\bf Remark.}~~The upper Lipschitz bound $\sqrt{B}$ has been obtained independently in \cite{BCMN13a}. The lower Lipschitz bound
we obtained here strenghtens the estimates given in \cite{BCMN13a}. 
Specifically their estimate for $\rho_{\infty}$ reads
$\sigmabf\leq\rho_{\infty}\leq\sqrt{2}\sigmabf$ where
\[ \sigmabf = \min_{S} \max(\sigma_n(F_S),\sigma_n(F_{S^c})) \]
Clearly $\sigmabf\leq \Delta\leq \sqrt{2}\sigmabf$. 

\vspace{3mm}

We conclude this section by turning our attention to the analysis of $V(x,y)$. A motivation
for studying it is that in practical problems the noise is
often added directly to $\nonlin^2$ as in (\ref{eq:noise1}) rather than to $\nonlin$. Such noise model is
used in many studies of phaseless reconstruction, e.g. in the Phaselift algorithm \cite{CSV12}, or in the IRLS algorithm in \cite{Bal12a}.

Let $\Sym_n(\R)$ denote the set of $n\times n$ symmetric matrices over $\R$. It is
a Hilbert space with the standard inner product given by $\ipT{X}{Y}:= \tr(XY^T)=\tr(XY)$.
The nonlinear map $\nonlin^2$ actually induces a linear map on $\Sym_n(\R)$. Write
$X=xx^T$ for any $x\in\R^n$. Then the entries of $\nonlin^2(x)$ are
\begin{equation} \label{4.13}
    (\nonlin^2(x))_j = |\inner{f_j,x}|^2=x^Tf_jf_j^T x = tr(F_jX) =\ipT{F_j}{X},
\end{equation}
where $F_j :=f_jf_j^T$. Now we denote by $\Ac$ the linear operator
$\Ac:~\Sym_n(\R)\longrightarrow\R^m$ with entries
$$
    (\Ac(X))_j = \ipT{F_j}{X} =\tr(F_j X).
$$
Let $S^{p,q}_n$ be the set of $n\times n$
real symmetric matrices that have at most $p$ positive and $q$
negative eigenvalues. Thus $S_n^{1,0}$ denotes the set of $n\times n$
real symmetric non-negative definite matrices of rank at most one.
Note that spectral decomposition easily shows that $X\in S_n^{1,0}$
if and only if $X=xx^T$ for some $x\in\R^n$.

The following lemma will be useful in this analysis
\begin{lemma}\label{lem-4.4}
The following are equivalent.
\begin{itemize}
\item[\rm (A)]~~$X\in S_n^{1,1}$.
\item[\rm (B)]~~$X = xx^T-yy^T$ for some $x,y\in\R^n$.
\item[\rm (C)]~~$X = \frac{1}{2}( w_1w_2^T + w_2w_1^T)$ for some $w_1, w_2 \in \R^n$.
\end{itemize}
Furthermore, for $X = \frac{1}{2}( w_1w_2^T + w_2w_1^T)$ its nuclear norm is
$\|X\|_1= \|w_1\|\|w_2\|$.
\end{lemma}
\Proof (A) $\Rightarrow$ (B) is a direct result of spectral decomposition, which yields
$X=\beta_1 u_1u_1^T - \beta_2 u_2u_2^T$ for some $u_1,u_2\in\R^n$ and
$\beta_1, \beta_2\geq 0$. Thus $X=xx^T-yy^T$ where $x:=\sqrt{\beta_1}u_1$
and $y:=\sqrt{\beta_2}u_2$.

(B) $\Rightarrow$ (C) is proved directly by setting $w_1=x-y$ and $w_2=x+y$.

We now prove (C) $\Rightarrow$ (A) by computing the eigenvalues of
$X  = \frac{1}{2}( w_1w_2^T + w_2w_1^T)$. Obviously $\rank(X) \leq 2$. Let $\lambda_1,\lambda_2$
be the two (possibly) nonzero eigenvalues of $X$. Then
$$
   \lambda_1+\lambda_2 = \tr\{X\} = \ip{w_1}{w_2},
$$
$$
    \lambda_1^2 + \lambda_2^2 = tr\{X^2\} = (\norm{w_1}^2\norm{w_2}^2+|\ip{w_1}{w_2}|^2)/2.
$$
Solving for eigenvalues we obtain
\begin{eqnarray*}
        \lambda_1 &=& \frac{1}{2}(\ip{w_1}{w_2} + \norm{w_1}\norm{w_2}),\\
        \lambda_2 &=& \frac{1}{2}(\ip{w_1}{w_2} - \norm{w_1}\norm{w_2}).
\end{eqnarray*}
Hence, by Cauchy-Schwartz inequality, $\lambda_1\geq 0 \geq \lambda_2$ which proves
$X\in S_n^{1,1}$. Furthermore, it also shows that the nuclear norm of $X$ is
$\|X\|_1 = |\lambda_1|+|\lambda_2| = \|w_1\|\|w_2\|$.
\eproof

Now we analyze $V(x,y)$. Parallel to the study of $U(x,y)$ we consider the following
quantities:
\begin{eqnarray*}
   \mu_\eps (x)  &:=&  \inf_{\{y: d(x,y)\leq \eps\}} V(x,y), \\
   \mu(x)  &:=&  \liminf_{\{y: d(x,y)\ra 0\}} V(x,y) = \liminf_{\eps\ra 0} \mu_\eps(x), \\
   \mu_0 &:=&  \inf_x \mu(x) , \\
   \mu_{\infty}  &:=&  \inf_{d(x,y)>0} V(x,y).
\end{eqnarray*}
as well as the upper bound $sup_{d_1(x,y)>0}V(x,y)$. 
By \eqref{4.13} we have $|\inner{f_j,x}|^2 -|\inner{f_j,y}|^2=\ipT{F_j}{X}$ where
$F_j = f_jf_j^T$ and $X=xx^T-yy^T$. Hence
$$
    V^2(x,y) = \frac{\sum_{j=1}^m |\ipT{F_j}{X}|^2}{\|X\|_1^2}.
$$
Set $w_1=x-y$ and $w_2=x+y$ and apply Lemma \ref{lem-4.4} we obtain
\begin{equation}  \label{4.14}
    V^2(x,y) =\frac{\sum_{j=1}^m |\ip{f_j}{w_1}|^2|\ip{f_j}{w_2}|^2}{\|w_1\|^2\|w_2\|^2}.
\end{equation}
We can immediately obtain the upper bound:
$$
   V(x,y) \leq \Bigl( \sup_{\norm{e_1}=1,\norm{e_2}=1} 
   \sum_{j=1}^m |\ip{f_j}{e_1}|^2 |\ip{f_j}{e_2}|^2 \Bigr)^{1/2}
 = \Bigl( \max_{\norm{e}=1} \sum_{j=1}^m |\ip{f_j}{e}|^4 \Bigr)^{1/2} =: \LTFnorm^2
$$
where $\LTFnorm$ denotes the operator norm of the linear analysis operator 
$T:H\rightarrow \R^m$, $T(x)=( \ip{x}{f_k} )_{k=1}^m $
defined between the Euclidian space $H=\R^n$ and the Banach space $\R^m$ endowed with 
the $l^4$-norm:
\begin{equation} \label{eq:r4}
\LTFnorm = \Bigl( \max_{\norm{x}=1} \sum_{k=1}^m |\ip{x}{f_k}|^4 \Bigr)^{1/4}
\end{equation}
Note also that
\[ \LTFnorm^2 = \max_{\norm{x}=1} \lambda_{max} (R(x)) \]
where $R(x)$ was defined in (\ref{eq:q2}). An immediate bound is $\LTFnorm \leq \sqrt{B}\max \norm{f_k}$
 with $B$ the upper frame bound of $\fc$.

Fix $x\neq 0$ and let $d(x,y) \ra 0$. Then either $y\rightarrow x$ or $y\ra -x$. Without loss
of generality we assume that $x\ra y$. Thus $w_1=x-y\rightarrow 0$ and $w_2=x+y\rightarrow 2x$.
However $w_1/\norm{w_1}$ can be any unit vector. Thus
$$
    \mu^2(x)  = \frac{1}{\norm{x}^2}\inf_{\norm{u}=1}\sum_{j=1}^m
    |\ip{f_j}{x}|^2 |\ip{f_j}{u}|^2 = \frac{1}{\norm{x}^2} \inf_{\norm{u}=1} \ip{R(x)u}{u}
    = \frac{1}{\norm{x}^2} \lambda_{\mathrm{min}}(R(x))
$$
where $R(x)$ was introduced in (\ref{eq:q2}).
Thus we obtain
$$
      \mu^2(x) = \frac{1}{\norm{x}^2}\lambda_{\mathrm{min}}(R(x)), \mhsp
       \mu^2_{0} = \min_{\norm{u}=1} \lambda_{\mathrm{min}}(R(u)).
$$
On the other hand note
$$
   \inf_{d(x,y)>0} V^2(x,y) = \inf_{w_1,w_2\neq 0}
    \frac{\sum_{j=1}^m |\ip{f_j}{w_1}|^2 |\ip{f_j}{w_2}|^2} {\norm{w_1}^2 \norm{w_2}^2 }
     = \min_{\norm{u}=1} \lambda_{\mathrm{min}}(R(u)) =a_0^2.
$$
where $a_0$ was introduced in (\ref{eq:a0}).
Thus we proved:

\begin{thm} \label{theo-4.5}
Assume the frame $\fc$ is phase retrievable. Then
\begin{eqnarray}
\mu(x) & = &  \frac{1}{\norm{x}} \sqrt{\lambda_{\mathrm{min}}(R(x))}, \\
\mu_{\infty} & = & \mu_0 =  \min_{u:\norm{u}=1} \sqrt{\lambda_{\mathrm{min}}(R(u))}=\sqrt{a_0}.
\end{eqnarray}
Furthermore $\nonlin^2$ is bi-Lipschitz with upper Lipschitz bound $\LTFnorm^2$ and lower Lipschitz bound $\mu_0$:
\[ \mu_0 d_1(x,y) \leq \norm{\nonlin^2(x)-\nonlin^2(y)} \leq \LTFnorm^2 d_1(x,y) \]
\end{thm}

\vspace{2mm}
\noindent
{\bf Remark.}~~ Note that the distance $d(.,.)$ is not equivalent to $d_1(.,.)$. Theorem
\ref{theo-4.5} now also implies that $\nonlin^2$ is not bi-Lipschitz with respect to
the distance $d(.,.)$ on $\hat H$. This fact was pointed out in \cite{BCMN13a}.

\section{Robustness and Size of Redundancy\label{sec5}}

\setcounter{equation}{0}

Previous sections establish results on the robustness of phaseless reconstruction
for the worst case scenario. A natural question is to ask: can ``reasonable''
robustness be achieved for a given frame, and in particular with small number of samples?
We shall examine how $q_\infty$ scales as the dimension $n$ increases.

Consider the case where $m=2n-1$. This is the minimal redundancy required for phaseless
reconstruction. In this case any frame $\fc$ would have $\Delta=\omega$. Hence we have
$\min\{1/\omega, 1/\eps\} \leq q_\eps = 1/\omega$. The stability of the reconstruction is
thus mostly controlled by the size of $1/\omega$. The question is: how big is $\omega$,
especially as $n$ increases?

Assume that the frame elements of $\fc$ are all bounded by $L$, $\|f_j\| \leq L$ for all
$f_j\in\fc$. Consider the $n+1$ elements $\{f_j:~j=1, \dots, n+1\}$. They are linearly dependent
so we can find $c_j\in\R$ such that $\sum_{j=1}^{n+1} c_j f_j =0$. Without loss of generality we may
assume $|c_{n+1}| = \min\{|c_j|\}$. Set $v = [c_1, c_2, \dots, c_n]^T$. Let $G =[f_1, \dots, f_n]$.
Then $Gv = \sum_{j=1}^n c_j f_j =-c_{n+1} f_{n+1}$. Now all $|c_j| \geq |c_{n+1}|$ so
$\|v\| \geq \sqrt n |c_{n+1}|$. Thus
$$
     \|Gv\| = |c_{n+1}| \|f_{n+1}\| \leq \frac{L}{\sqrt n} \|v\|.
$$
It follows that $\sigma_n(G) \leq \frac{L}{\sqrt n}$, and hence
\begin{equation}  \label{5.1}
    \omega \leq \frac{L}{\sqrt n}.
\end{equation}

Note that here we have considered only the first $n+1$ vectors of the frame $\fc$. The actual
value of $\omega$ will likely decay much faster as $n$ increases. In a preliminary 
work we are able
to establish the bound $\omega \leq CL/\sqrt{n^3}$ where $C$ is independent of 
$n$ \cite{Wan13}. But even this estimate is likely far from optimal.

\smallskip
\begin{conjecture}  \label{conj-5.1}
 Let $m =2n-1$ and $\|f_j\| \leq L$ for all $f_j\in\fc$. Then
 there exist constants $C>0$ and $0<\beta<1$ independent of $n$ such that
$$
     \omega \leq CL \beta^{-n}.
$$
\end{conjecture}

\smallskip

A related problem is as follows: Consider an $n\times (n+k)$ matrix 
$F=[g_1,g_2, \dots, g_{n+k}]$. Let $\tau = \min \{\sigma_n (F_S):~S\subset \{1,\dots, n+k\},
|S|=n\}$. Assume that all $\|g_j\| \leq 1$. How large can $\tau$ be? For $k=1$ we have
already seen that it is bounded from above by $C/\sqrt n$. The preliminary work \cite{Wan13}
shows that for $k=1$ it is bounded from above by $C/n^{\frac{3}{2}}$.

\begin{conjecture}  \label{conj-5.2}
 There exists a constant $C=C(k)$ such that
 $$
     \tau \leq \frac{C}{n^{k-\frac{1}{2}}}.
$$
\end{conjecture}

Thus in the minimal setting with $m=2n-1$ it is impossible to achieve scale independent stability for phaseless reconstruction. The same arguments can be used to show that
even when $m=2n+k_0$ for some fixed $k_0$ scale independent stability is not possible.
A natural question is whether scale independent stability is possible when we increase
the redundancy of the frame. As it turns out this is possible via a recent work
by Wang \cite{Wan12}. More precisely, the following result follows from the main results 
in \cite{Wan12}:

\begin{thm} \label{theo-5.3}
    Let $r_0 >2$ and let $F=\frac{1}{\sqrt n}G$ where $G$ is 
    an $n\times m$ random matrix whose elements are i.i.d.
    normal $N(0,1)$ random variables such that $m/n = r_0$. Then there exist
    constants $0<\Delta_0 \leq \omega_0$ dependent only on $r_0$ and not on $n$ such that
    with high probability we have
    $$
         \Delta \geq \Delta_0, \mhsp  \omega \geq \omega_0.
    $$
\end{thm}
\Proof  Part of the main theorem of \cite{Wan12} proves the following result: 
Let $\lambda >\delta>1$ be fixed. Assume that 
$A=\frac{1}{\sqrt{n}}B$ where $B$ is an $n \times N$ random Gaussian matrix
with i.i.d $N(0,1)$ entries such that $N/n = \lambda$. Then there exists a constant
$c>0$ depending continuously and only on $\lambda$ and $\delta$ such that
$$
   \min_{S\subseteq\{1,\dots,N\}, |S|\geq\delta n} \sigma_n(A_S) \geq c.
$$
The theorem now readily follows. Observe that because $r_0>2$, in the expression for
$\Delta$ we may choose $\lambda= r_0$ $\delta = \frac{r_0}{2}>1$ and clearly we have
$$
   \Delta \geq \min_{S\subseteq\{1,\dots,N\}, |S|\geq\delta n} \sigma_n(F_S) \geq \Delta_0.
$$
for some $\Delta_0>0$ independent of $n$. For $\omega$ we may choose $\lambda=r_0$ 
and $\delta = r_0-1>1$. Again the theorem of \cite{Wan12} implies that
$$
   \omega \geq \min_{S\subseteq\{1,\dots,N\}, |S|\geq\delta n} \sigma_n(F_S) \geq \omega_0.
$$
\eproof

\smallskip

In the theorem the values $\Delta_0$ and $\omega_0$ can be estimated explicitly. Here with
high probability is in the standard sense that the probability is at least 
$1- c_0 e^{-\beta n}$ for some $c_0,\beta>0$. Thus scale independent stable phaseless
reconstruction is
possible whenever the redundancy is greater than $2+\delta$, $\delta>0$, at least for
random Gaussian matrices.

\vspace{4mm}
\noindent
{\bf Acknowledgement.} ~~The authors would like to thank Matt Fickus, Dustin Mixon and
Jeffrey Schenker for very helpful discussions.

\end{document}